# A Note on Thompson Problem

In celebration of the 90 birthday of Professor John Thompson

Yu Li [a] and Wujie Shi [b],*
[a] School of Mathematics and Statistics, Guangxi Nornal University, Guilin 541004, China
School of Mathematics and Statistics, Southwest University, Chongqing 400715, China
and
[b] Chongqing University of Arts and Sciences, Chongqing 402160, China;
Suzhou University, Jiangsu 215006, China
e-mails: liskyu@163.com and shiwujie@outlook.com

**Abstract**: In this paper, we prove that if two finite groups $G$ and $H$ have isomorphic Burnside rings, then $G$ and $H$ are the same order type groups, and give an example to show that the Burnside rings of the same order type groups are not necessarily isomorphic. This result is related to the Thompson Problem which was raised in 1987.

# 1  Introduction

In 1987, the second author conjectured that every finite simple group is recognizable in the class of finite groups by its order and the set of element orders. That is:

**Conjecture 1**. Let $G$ be a finite group and $S$ a finite simple group. Then $G \cong S$ if and only if (1) $\pi_e(G) = \pi_e(S)$, where $\pi_e(G)$ denote the set of element orders of $G$, and (2) $|G| = |S|$.

After reported the above conjecture to Professor J.G. Thompson, Thompson posed the following problem in his reply [17]:

For each finite group $G$ and each integer $d > 1$, let $G(d) = \{x \in G \mid x^d = 1\}$.

**Definition 2.** $G_1$ and $G_2$ are of the same order type if and only if $|G_1(d)| = |G_2(d)|$, $d = 1, 2, \ldots$.

**Thompson Problem (1987).** Suppose $G_1$ and $G_2$ are groups of the same order type. Suppose also that $G_1$ is solvable. Is it true that $G_2$ is also necessarily solvable?

---

\* Corresponding author.

Moreover, Thompson pointed out that "**I have talked with several mathematicians concerning groups of the same order type. The problem arose initially in the study of algebraic number fields, and is of considerable interest.**"

In fact, the groups of the same order type are such groups whose order equations are the same. That is, regarding the element orders as the equivalent relation, the group $G$ is partitioned into the same order classes. The order equation is

$$|G| = n_1 + n_2 + \ldots + n_s,$$

where $n_i > 0$ denotes the number of elements of order $i$ in $G$.

For convenience, $|G|$ and $\pi_e(G)$ are said to be **two orders** of $G$. After the publications of a series of papers (see [4, 10, 11, 12, 13, 14, 15, 18, 19]), Conjecture 1 had proved and become a theorem. That is, all finite simple groups $G$ can be determined by their two orders. Moreover, symmetric groups are all characterized by their two orders[2], but dihedral group $D_8$ and quaternion group $Q_8$ are not characterized by their two orders. So, we have the following question:

**Question 3**. What kind of finite groups can be characterized by their two orders?

**Definition 4**. $G_1$ and $G_2$ are of *the same two orders type* if and only if $|G_1| = |G_2|$ and $\pi_e(G_1) = \pi_e(G_2)$.

Up to now we do not answer Thompson Problem. Considering the counter-example, the same order type groups which are not isomorphic, Thompson has also given a non-solvable example of same order type groups which are not isomorphic: $G_1 = 2^4 : A_7$ and $G_2 = L_3(4) : 2_2$ (the notation from ATLAS [5]) both are maximal subgroups of $M_{23}$.

Comparing Thompson problem we put forward the following question: Suppose $G_1$ and $G_2$ are the same two orders type. Suppose again that $G_1$ is solvable. Is it true that $G_2$ is also necessarily solvable?

There is a counter-example for the above question.

**Examples 5.** Let $G_1 = A_5 \times A_5 \times A_5$, $G_1$ is not solvable and $|G_1| = 60^3$ and $\pi_e(G_1) = \{1, 2, 3, 5, 6, 10, 15, 30\}$ since $\pi_e(A_5) = \{1, 2, 3, 5\}$. Let $G_2 = C_{30} \times C_{30} \times C_{30} \times C_2 \times C_2 \times C_2$. Then $|G_2| = 60^3$ and $\pi_e(G_2) = \{1, 2, 3, 5, 6, 10, 15, 30\}$. $G_1$ and $G_2$ are the same two orders type, but $G_2$ is solvable. Moreover, $G_1$ and $G_2$ are not the same order type.

Let $G$ be a finite group and $B(G)$ the Burnside ring of $G$ as defined [1, 6]. Dress [7] proved that if two finite groups $G$ and $H$ have isomorphic Burnside rings, then $|G| = |H|$. It is natural to ask what properties does two groups have if there Burnside rings are isomorphic. It was asked in [20] whether two groups are isomorphic if their Burnside rings are isomorphic, Thevenaz [16] obtained the first counter-example of this problem. In [9], it was shown that if two finite groups have isomorphic Burnside rings, then there is a one-to-one correspondence between their families of soluble subgroups which preserves order and conjugacy class of subgroups, Kimmerle, Luca and Raggi-Cardenas [8] proved that if two finite groups have isomorphic Burnside rings, then there is a one-to-one correspondence between their families of cyclic subgroups which preserves order and conjugacy class of subgroups. It turns out that if two finite groups $G$ and $H$ with isomorphic Burnside ring, then $\pi_e(G) = \pi_e(H)$ [8, Corollary 5.2]. If $G$ is a finite simple group and $B(G) \cong B(H)$, then by Conjecture 1, $G \cong H$. Therefore, isomorphism class of the simple group is completely determined by its Burnside ring, and the answer of Yoshida's problem [20, problem 2] is affirmative for simple groups. In this paper, we prove that if two finite groups $G$ and $H$ have isomorphic Burnside rings, then $G$ and $H$ are the same order type, and give an example to show that the Burnside rings of the same order type groups are not necessarily isomorphic.

## 2 The Results

The following two lemmas are important for the proof our main results.

**Lemma 6.** [*9, Corollary 5.3*] *Let $G$ and $G'$ be finite groups such that their Burnside rings are isomorphic. There is one-to-one correspondence between the conjugacy classes of solvable subgroups of $G$ and $G'$ which preserves order of subgroup and cardinality of the conjugacy class.*

**Lemma 7.** [*8, Proposition 5.1*] *Let $G$ and $H$ be finite groups and $a : B(G) \to B(H)$ be an isomorphism between their Burnside rings. Then there is an order- and inclusion- preserving bijection $\sigma_*$ between the conjugacy classes of nilpotent subgroups of $G$ and $H$. Moreover, $\sigma_*$ maps a class of cyclic subgroups to a class consisting of cyclic subgroups.*

We are now ready to prove main results.

**Theorem 8.** Let $G$ and $H$ be finite groups such that their Burnside rings are isomorphic. Then $G$ and $H$ are the same order type.

*Proof.* If the Burnside rings of $G$ and $H$ are isomorphic, by Lemma 7, there is a bijection $\sigma_*$ between the conjugacy classes of cyclic subgroups of $G$ and $H$ which preserves order and cardinality of the conjugacy class. If $N_1, N_2, \ldots, N_s$ and $K_1, K_2, \ldots, K_t$ are the representation of conjugacy classes of cyclic groups of order $n$ of $G$ and $H$, respectively, then $s = t$ and $|G : N_G(N_i)| = |H : N_H(K_j)|$ if $\sigma_*(N_i) = K_j$, $i, j = 1, 2, \ldots, s$. It follows that $\Sigma|G : N_G(N_i)|\varphi(n) = \Sigma|H : N_H(K_i)|\varphi(n)$, where $\varphi$ is Euler's function. It implies that $G$ and $H$ have the same number of elements of order $n$. Therefore, $G$ and $H$ are the same order type. □

In the following theorem, we show that the converse of Theorem 8 is not true.

**Theorem 9.** *Let $G_1 = L_3(4) : 2_2$ and $G_2 = 2^4 : A_7$. Then $G_1$ and $G_2$ are the same order type*, but $B(G_1)$ *is not isomorphic to* $B(G_2)$.

*Proof.* Computation in MAGMA [3] shows that $G_1$ and $G_2$ have the same order equation

$$|G_1| = |G_2| = n_1 + n_2 + n_3 + n_4 + n_5 + n_6 + n_7 + n_8 + n_{14}$$
$$= 1 + 435 + 2240 + 6300 + 8064 + 6720 + 5760 + 5040 + 5760.$$

So $G_1$ and $G_2$ are the same order type. By ATLAS [5], $A_7$ has a solvable maximal subgroup $(A_4 \times 3) : 2$ of order $2^3 \cdot 3^2$ and $2^3 \cdot 3^2$ doesn't divide the order of all maximal subgroups of $L_3(4)$. It follows that $G_2$ has a solvable subgroup of order $2^7 \cdot 3^2$, but $G_1$ has no subgroup of order $2^7 \cdot 3^2$. Hence by Lemma 5, $B(G_1)$ *is not isomorphic to* $B(G_2)$. □

In this note, we investigated the following three properties of two finite groups $G_1$ and $G_2$: (1) $B(G_1) \cong B(G_2)$, that is, $G_1$ and $G_2$ are the same Brunside ring type, (2) $G_1$ and $G_2$ are the same order type, (3) $G_1$ and $G_2$ are are the same two orders type. From the above discussion, we get that these three statements are not equivalent to each other, but (1) implies (2) and (2) implies (3). The Thompson Problem is that two finite groups satisfy property (2), one solvable, and is the other also solvable? If (2) is changed to (1), the answer is affirmative by [6], while (2) is changed to (3), the answer is negative and Example 5 is a counter-example.

# Acknowledgements

This research was supported by Guangxi Natural Science Foundation (2022GXNSFBA035572) and Natural Science Foundation of China (11171364, 11271301).

# References

[1] D.J. Benson, Representations and Cohomology I, in: Cambridge Stud. Adv. Math., vol. 30, Cambridge Univ. Press, Cambridge, 1991.

[2] J.X. Bi, A new characterization of the symmetry groups. *Acta Math. Sinica* **33** (1) (1990), 70-77(in Chinese).

[3] W. Bosma, J. Cannon, C. Playoust, The Magma algebra system I. The user language, *J. Symbolic Computation* **24** (1997), 235-265.

[4] H.P. Cao and W.J. Shi, Pure quantitative characterization of finite projective special unitary groups, *Sci. China, Ser.* A **45** (6) (2002), 761-772.

[5] J.H. Conway, R.T. Curtis, S.P. Norton, R.A. Parker, R.A. Wilson, Atlas of Finite Groups, Oxford University Press, London, 1985.

[6] A. Dress, A characterization of solvable groups, *Math. Z.* **110** (1969), 213-217.

[7] A. Dress, Notes on the theory of representations of finite groups. Part I: The Burnside ring of a finite group and some AGN-applications. With the aid of lecture notes, University of Bielefeld, 1971.

[8] W. Kimmerle, F. Luca and A.G. Raggi-Cardenas, Irreducible components and isomorphisms of the Burnside ring, *J. Group Theory* **11**(2008), 831-844.

[9] A.G. Raggi-Cardenas and L. Valero-Elizondo, Normalizing isomorphisms between Burnside rings, *J. Algebra* **277** (2004), 643-657.

[10] W.J. Shi, A new characterization of the sporadic simple groups. Group Theory: Proceedings of the 1987 Singapore Conference, Walter de Gruyter, Berlin, 531-540 (1989).

[11] W.J. Shi, The pure quantitative characterization of finite simple groups (I), *Progr. Natur. Sci.* **4** (3) (1994), 316-326.

[12] W.J. Shi and J.X. Bi, A characteristic property for each finite projective special linear group, Lecture Notes in Math., Springer-Verlag, **1456** (1990), 171-180.

[13] W. J. Shi and J.X. Bi, A characterization of Suzuki-Ree groups, *Sci. in China, Ser. A* **34** (1991), 14-19.

[14] W.J. Shi and J.X. Bi, A new characterization of the alternating groups, *Southeast Asian Bull. Math.*, **16** (1) (1992), 81-90.

[15] R. Solomon, On finite simple groups and their classification, *Notices Amer. Math. Soc.* **42** (2) (1995), 231-239.

[16] J. Thaevenaz. Isomorphic burnside rings, *Comm. Alg.* **16** (1988), 1945-1947.

[17] J.G. Thompson, Private communication, 1987.

[18] A.V. Vasil'ev, M.A. Grechkoseeva, and V.D. Mazurov, Characterization of the finite simple groups by spectrum and order, *Algebra Logic* **48** (6) (2009), 385-409.

[19] M.C. Xu and W.J. Shi, Pure quantitative characterization of finite simple groups $^2D_n(q)$ and $D_l(q)$ ($l$ odd), *Algebra Colloq.* **10** (3) (2003), 427-443.

[20] T. Yoshida, On the Burnside rings of finite groups and finite categories. Commutative algebra and combinatorics (Kyoto, 1985), 337-353, *Adv. Stud. Pure Math.*, **11** (North-Holland, 1987), 337-353.